\documentclass[12pt,draft]{article}

\usepackage{bm}        
\usepackage{amsmath}   
\usepackage{amssymb}   
\usepackage{stmaryrd}  
\usepackage{eucal}     

\newtheorem{lemma}{Lemma}
\newtheorem{theorem}{Theorem}
\newtheorem{corollary}{Corollary}
\newtheorem{pro}{Proposition}
\newtheorem{theoremL}{Theorem}

\def \s {\sigma}

\def \r {\rho}
\def \la {\langle}
\def \ra {\rangle}
\def \M {\mathcal{M}}
\def \D {\mathcal{D}}
\def \E {\mathcal{E}}
\def \I {\mathcal{I}}

\def \a {\alpha}

\def \lm {\lambda}
\def \ov {\overline}

\def \ll {\llbracket}
\def \rr {\rrbracket}
\def \ZZ {\mathbb{Z}}
\def \be {\begin{equation}}
\def \ee {\end{equation}}
\def \ve {\varepsilon}
\def \F {\bm{F}}
\newcommand{\ba}[1] {\begin{array}{#1}}
\newcommand{\ea} {\end{array}}
\def \llb {\llbracket}
\def \rrb {\rrbracket}
\def \O {\mathbb{O}}
\def \ph {\phantom}

\begin{document}

\begin{center}
{\Large \bf Sylow's theorem for Moufang loops}

\vspace{0.5cm}

{\sc Alexander N.\,Grishkov\footnote{Departamento de Matem\'atica,
USP, S\~ao Paulo, Brazil.}, \ Andrei
V.\,Zavarnitsine\footnote{Sobolev Institute of Mathematics,
Novosibirsk, Russia}}
\end{center}

\begin{abstract}
For finite Moufang loops, we prove an analog of the first Sylow theorem giving a criterion of the
existence of a $p$-Sylow subloop. We also find the maximal order of $p$-subloops in the Moufang loops
that do not possess $p$-Sylow subloops.
\end{abstract}

\begin{center}{\small
}
\end{center}

\section{Introduction}

After a positive solution to the Lagrange problem for finite
Moufang loops was given \cite{gz_lagr}, proving an analog of
Sylow's theorems has become an important problem in the theory of
finite Moufang loops \cite{gz_syl_pr}. We should mention right
away that an obvious obstruction to the precise analog of the
first Sylow theorem about the existence of a $p$-Sylow subloop has
been known for a long time: if $M(q)$ is a finite simple Paige
loop over a finite field $\F_q$ then the orders of elements of
$M(q)$ divide $q(q^2-1)$, whereas the order of $M(q)$ is
$\frac{1}{d}q^3(q^4-1)$, where $d=\gcd(2,q-1)$. Hence, for a prime
$p$ dividing $\frac{1}{d}(q^2+1)$, a $p$-Sylow subloop of $M(q)$
does not exist; consequently, nor does it exist in any finite
Moufang loop having a composition factor isomorphic to $M(q)$.

The main result of the present paper asserts that this obstruction
to the existence of a $p$-Sylow subloop is in fact the only one.
Namely, if a composition series of a Moufang loop $M$ contains no
simple factors $M(q)$ such that $p\,|\,\frac{1}{d}(q^2+1)$, then
$M$ has a $p$-Sylow subloop. In this case, we call $p$ a {\it
Sylow prime} for $M$. In particular, every prime is Sylow for all
solvable finite Moufang loops, i. e., a $p$-Sylow subloop exists
in this case for every $p$. We note that this result for Moufang
loops of odd order (in fact, a stronger result on the existence of
Hall $\pi$-subloops) was proved in \cite{gla2}.

Many major results about finite Moufang loops were proved based on
the correspondence between Moufang loops and groups with triality
\cite{dor}. This paper is not an exception. Clearly, the existence
of a $p$-Sylow subloop is related to the existence of an
$S$-invariant subgroup in a corresponding group $G$ with triality
$S$. Unfortunately, even when the above obstruction to the
existence of a $p$-Sylow subloop in $M$ is absent (which is always
the case, for example, when $p=2$ or $3$), one cannot guarantee
existence in $G$ of an $S$-invariant $p$-Sylow subgroup. We
therefore prove only that $G$ possesses a sufficiently large
$S$-invariant $p$-subgroup whose corresponding $p$-subloop in $M$
is $p$-Sylow.

A Moufang loop $M$ may still have nontrivial $p$-subloops even
when $p$ is not a Sylow prime for $M$. The order of such a
$p$-subloop is obviously bounded by the product of the orders of
$p$-Sylow subloops of the composition factor of $M$ for which $p$
is Sylow. In the last section, we show that this bound is actually
achieved for all finite Moufang loops $M$. We call such
$p$-subloops of maximal order {\it quasi-$p$-Sylow}. In this way,
Sylow's theorem can be reformulated as follows:

\begin{theoremL}\label{tha} Every finite Moufang loop has a quasi-$p$-Sylow
subloop for all primes~$p$.
\end{theoremL}
The proof of this result uses the following important structural
theorem (see Theorem \ref{MT} in section 5) for Moufang loops:

\begin{theoremL}\label{thb} Every finite Moufang loop $M$ contains uniquely determined normal
subloops $Gr(M)$ and $M_0$ such that $Gr(M)\leqslant M_0$, $M/M_0$ is an elementary abelian $2$-group,
$M_0/Gr(M)$ is the direct product of simple Paige loops $M(q)$ (where $q$ may vary), the composition
factors of $Gr(M)$ are groups, and $Gr(M/Gr(M))=1$.
\end{theoremL}

Another interesting auxiliary fact we obtain is that a
nonassociative minimal normal subloop of a finite Moufang loop
must necessarily be simple (see (i) of Theorem \ref{MT}). This
extends the well-known group-theoretic assertion that a minimal
normal subgroup of a finite group is the direct product of
isomorphic simple groups.

It is natural to conjecture that the analogs of other two Sylow
theorems, namely the embedding of $p$-subloops of $M$ into a
quasi-$p$-Sylow subloop and the conjugacy of quasi-$p$-Sylow
subloops by (inner) automorphisms of $M$, are true for finite
Moufang loops as well. As far as determining the number of
quasi-$p$-Sylow subloops, at present we do not have a formulation
of the corresponding conjecture for lack of sufficient
experimental data.

\section{Preliminaries}

All loops we consider are finite. If $x,y$ are elements of a group
$G$ then $x^y=y^{-1}xy$, $x^{-y}=(x^{-1})^y$, $[x,y]=x^{-1}x^y$.
$Z(G)$ is the center of $G$. If $\varphi\in Aut(G)$ then
$x^\varphi$ is the image of $x$ under $\varphi$. If $n$ is an
integer and $p$ is a prime then the notation $p^k\, \|\, n$ for
$k\geqslant 0$ means that $p^k \mid n$ and $p^{k+1} \nmid n$. For
natural numbers $m,n$, we write $(m,n)$ for $\gcd(m,n)$.

A loop $M$ is called a {\it Moufang loop} if
$xy\cdot zx = (x\cdot yz)x$ for all $x,y,z \in M$.
For $x,y,z\in M$ define the {\it commutator} $\ll x,y \rr$
by $xy=yx\cdot \ll x,y \rr$.
The {\it nucleus}
$Nuc(M)$ of a Moufang loop $M$ is the set
$\{a\in M \ | \ a\cdot xy=ax \cdot y\ \ \forall x,y\in M\}$.
For basic properties of Moufang loops, see \cite{bru}.

To arbitrary elements $x,y$ of a Moufang loop $M$, there are associated bijections $L_x$, $R_x$, $T_x$,
$L_{x,y}$, $R_{x,y}$ of $M$ defined as follows:

\be \ba{c}
 yL_x=xy, \quad yR_x=yx, \quad \text{for all} \quad y\in M;\\
 T_x=L_x^{-1}R_x, \qquad
L_{x,y}=L_xL_yL_{yx}^{-1}, \qquad R_{x,y}=R_xR_yR_{xy}^{-1}.
\ea \ee
The {\it multiplication group} $Mlt(M)$ of $M$ is the group
of permutations of $M$ generated by $L_x$ and
$R_x$ for all $x \in M$, and the {\it inner mapping group}
$\I(M)$ is the subgroup of $Mlt(M)$ generated by $T_x$ and
$R_{x,y}$ for all $x,y$  in $M$.

A {\it pseudoautomorphism} of a Moufang loop
$M$  is a bijection $A: M\to M$
with the property that there exists an element $a
\in M$ such that
$$xA(yA\cdot a)=(x\cdot y)A \cdot a \qquad \text{for all} \quad x,y \in M.$$
Such an element $a$ is called a {\it right companion} of $A$.
Denote by $PsAut(M)$ the group formed by
all pairs $(A,a)$, where $A$ is a pseudoautomorphism of $M$ with right
companion $a$, with respect to the operation
$(A,a)(B,b)=(AB,aB\cdot b)$.
Denote by $PsInn(M)$ the
subgroup of $PsAut(M)$ generated by the elements
$(T_x,x^{-3})$ and $(R_{x,y},\llb x,y \rrb)$ for all $x,y$ in $M$
(such elements are in $PsAut(M)$ by Lemma VII.2.2 in \cite{bru}).

By \cite{lie}, the only finite simple nonassociative Moufang loops are the Paige loops $M(q)$ which
exist for every finite field $\F_q$. The order of $M(q)$ is $\frac{1}{d}q^3(q^4-1)$, where $d=(2,q-1)$,
and is the product of two coprime numbers $q^3(q^2-1)$ and $\frac{1}{d}(q^2+1)$.

Let $M$ be a Moufang loop and let $p$ be a prime. $M$ is a {\it $p$-loop} if the order of every element
of $M$ is a power of $p$ (for Moufang loops, this is equivalent to the condition that $|M|$ be a power
of $p$). A {\it $p$-Sylow subloop} of $M$ is a subloop of order $p^k$, where $p^k\, \|\ |M|$,
$k\geqslant 0$. We denote by $Syl_p(M)$ the set of all $p$-Sylow subloops of $M$. The fact that
$Syl_p(M)$ can be empty for $p$ dividing $|M|$ was observed long ago. For example, the simple Paige loop
$M(2)$ of order $120$ does not have elements of order $5$. The classification \cite{gz_max} of maximal
subloops of $M(q)$ implies the following assertion:

\begin{lemma} \label{syl_simple} If $p\nmid \frac{1}{d}(q^2+1)$ then $Syl_p(M(q))\ne \varnothing$.
If $p\ |\ \frac{1}{d}(q^2+1)$ then $M(q)$ does not have elements of order $p$; in particular,
$Syl_p(M(q))=\varnothing$ in this case.
\end{lemma}

{\it Proof.} First, suppose that $q$ is a power of $p$. Consider a maximal subloop of $M(q)$ of shape
$q^2:PSL_2(q)$ (i.e. a split extension of a normal elementary abelian subgroup order $q^2$ by the simple
group $PSL_2(q)$). This subloop obviously contains a subloop of order $q^3$ which is therefore
$p$-Sylow in $M(q)$.

If $p\ |\ q^2-1$ then consider a maximal subloop of type $M(PSL_2(q),2)$ (i.e. a loop with a
normal subgroup $PSL_2(q)$ of index $2$ constructed by the Chein duplication process as in
Theorem 1 of \cite{che}). Then obviously $Syl_p(PSL_2(q))\subseteq Syl_p(M(q))$ whenever $p$
is odd, and $M(S,2)\in Syl_2(M(q))$ for every $S\in Syl_2(PSL_2(q))$.

If $p\ |\ \frac{1}{d}(q^2+1)$ then it can be easily seen by induction that every maximal subloop of
$M(q)$ either has order coprime with $p$ or contains no elements of order $p$. $\blacktriangle$

Note that we will give another proof (see Corollary \ref{alt} below) of the existence of a
$p$-Sylow subloop in $M(q)$ which does not use the classification  \cite{gz_max}. However, the
explicit structure of the Sylow subloops of $M(q)$ is best seen from the proof of Lemma
\ref{syl_simple} above.

Lemma \ref{syl_simple} also shows that there is an obstruction to the existence of a $p$-Sylow subloop
in an arbitrary finite Moufang loop $M$ having a composition factor $M(q)$ with $\frac{1}{d}(q^2+1)$
divisible by $p$.

{\it Definition.} Let $M$ be a Moufang loop. A prime $p$ is called a {\it Sylow prime for $M$} if, for
every composition factor of $M$ that is isomorphic to $M(q)$ for some $q$, we have $p\ \nmid
\frac{q^2+1}{(2,q-1)}$.

Obviously, $p$ is a Sylow prime for $M$ if and only if it is
such for all composition
factors of $M$. We can now state the main theorem.

\begin{theorem}[Sylow's theorem] \label{main}
Let $M$ be a finite Moufang loop and let $p$ be a prime. Then $M$ contains
a $p$-Sylow subloop if and only if $p$ is a Sylow prime for $M$.
\end{theorem}

Before proving this theorem, we will require some facts about groups with triality which are introduced
in the next section. As an important corollary to Sylow's theorem, we obtain the following assertion:

\begin{corollary} \label{cor1}
(i)\ Every Moufang loop has a $2$-Sylow and a $3$-Sylow subloop.

(ii)\  If all composition factors of a Moufang loop $M$ are groups then $M$ has a $p$-Sylow subloop for
all primes $p$.
\end{corollary}

{\it Proof.} (i)\ The primes $2$ and $3$ are Sylow for every loop, since $\frac{q^2+1}{(2,q-1)}$ is
coprime with $6$ for all prime powers $q$.

(ii)\ If $M$ is as stated then all primes are Sylow for $M$ by definition. $\blacktriangle$

\section{Groups with triality}

A group $G$ possessing automorphisms $\r$ and $\s$ that
satisfy $\r^3=\s^2=(\r\s)^2=1$ is called a {\it group
with triality $\la\r,\s\ra$} if
\begin{equation} \label{trial}[x,\s][x,\s]^\r[x,\s]^{\r^2}=1 \end{equation}
for every $x$ in $G$, where $[x,\s]=x^{-1}x^{\s}$. We henceforth
denote $S=\la \s,\r \ra$. Obviously, $S$ is a homomorphic image of
the symmetric group $S_3$ of degree $3$. \big(Strictly speaking,
we have fixed a homomorphism $\gamma:S_3\to Aut(G)$ and then
identified the generators $\r$ and $\s$ of $S_3$ with their images
under $\gamma$ to simplify the notation. We thus implicitly
consider $G$ as a {\it group with operators} $S_3$ in the sense of
\S15 of \cite{kur}.\big) It should be noted that the identity
(\ref{trial}) does not depend on a particular choice of the
generators $\r$ and $\s$ \cite{dor}.

Let $G$ be a group with triality $S$. Put
$M=\{ [x,\s]\ |\ x\in G\}$ and $H=C_G(\s)$.
Then $M$ endowed with the multiplication
\begin{equation} \label{loop_mult}
m.n=m^{-\r} n  m^{-\r^2} \ \ \text{for all} \ \ m,n\in M
\end{equation}
becomes a Moufang loop $(M,.)$ of order $|G:H|$. We denote this loop by $\M(G)$. Every Moufang
loop can be obtained in this way from a suitable group with triality. Moreover, for every
subloop $M_0\leqslant \M(G)$, there exists an $S$-invariant subgroup $G_0$ of $G$ (in brief,
{\it $S$-subgroup}) such that $M_0=\M(G_0)$. Any element of $\M(G)$ has the same order whether
viewed as a group or loop element. For more details on this, see
\cite{gz_lagr,gz_tri}. We observe that
\be \label{m-1}
m^{\s}=m^{-1}\in M \ \ \ \text{for all} \ \ \ m\in M .
\ee

A homomorphism $\varphi: G_1\to G_2$ of groups $G_1$ and $G_2$
with triality $S$ is called an {\it $S$-homo\-morphism} if
$\a\varphi=\varphi\a$ for all $\a \in S$. \big(Again, strictly
speaking, we have fixed homomorphisms $\gamma_i: S_3\to Aut(G_i)$,
$i=1,2$, and then required
$(\a\gamma_1)\varphi=\varphi(\a\gamma_2)$ for all $\a \in S_3$
thus making $\varphi$ an {\it operator homomorphism} from $G_1$ to
$G_2$. This is clearly a morphism in the category of groups with
triality.\big) Denote by $Z_S(G)$ the {\it $S$-center} of $G$,
which is by definition the maximal normal $S$-subgroup of $G$ on
which $S$ acts trivially.

\begin{lemma}\label{dual} Let $G$ be a group with triality $S$ and let $M=\M(G)$.
Then

(i) $[[G,S],S]=[G,S]$,

(ii) $Z_S(G/Z_S(G))=1$ and $Z_S(G)=C_G([G,S]S)$,

(iii) $[G,S]$ is generated by $M\cup M^\r$,

(iv) the elements $m$, $m^\r$, $m^{\r^2}$ of $G$
pairwise commute for all $m\in M$,

(v)  $m^{-\r} n  m^{-\r^2}=n^{-\r^2} m  n^{-\r}$ \ \ for all $m,n\in M$,

\end{lemma}

{\it Proof.} See Lemmas 1 and 2 in \cite{gz_tri}.  $\blacktriangle$

\begin{lemma}\label{tripl} Let $G$ be a group with triality $S=\la \s,\r \ra$
and let $M=\M(G)$. Then the semidirect product
$G_0=G\rtimes \la \r \ra$ is a group with triality $S$ if and only if $M$
has exponent~$3$.

\end{lemma}
{\it Proof.} For all $g\in G$, we have
$(g\r)^{-1}(g\r)^{\s}=\r^2g^{-1}g^{\s}\r^2$ and
$(g\r^2)^{-1}(g\r^2)^{\s}=\r g^{-1}g^{\s}\r$. Put $m=g^{-1}g^\s$.
Then
$$(\r m\r)(\r m\r)^\r(\r m\r)^{\r^2}=m^{\r^2}m^\r m=1,$$
$$(\r^2 m\r^2)(\r^2 m\r^2)^\r(\r^2 m\r^2)^{\r^2}=(m^\r)^3.$$
Hence, $G_0$ is a group with triality iff $m^3=1$ for all $m\in \M(G)$. $\blacktriangle$

Note that if $\M(G)$ has exponent $3$ then $G\rtimes S$ is also a group with triality $S$ (where $S$
acts by inner automorphisms). However, in this case $\M(GS)=\M(G\la\r\ra)$.

Given a Moufang loop $M$, there exists a universal group with
triality $\D(M)$ introduced by Doro \cite{dor} which satisfies
$\M(\D(M))\cong M$, $[\D(M),S]=\D(M)$ and, if $G$ is any group with
triality such that $\M(G)\cong M$ and  $G=[G,S]$,
then there is an $S$-epimorphism $\tau:\D(M)\to G$.
Denote $\E(M)\cong \D(M)/Z_S(\D(M))$. Then $\M(\E(M))\cong M$.

A group $G$ (with triality) is called {\it $S$-simple} if it has
no proper $S$-homo\-morphic images or, equivalently, contains no
normal $S$-subgroups.

\begin{lemma} \label{s-sim} Every finite $S$-simple group
is $S$-isomorphic to one of the following $S$-simple groups:

(i) A finite simple group $G$ with trivial $S$-action. In this case, $\M(G)$
is the trivial loop;

(ii) $G=\la a\ |\ a^3=1\ra\cong \ZZ_3$, $a^\r=a$, $a^\s=a^{-1}$.
In this case, $\M(G)\cong \ZZ_3$;

(iii) $G=\la a,b\ |\ a^p=b^p=[a,b]=1 \ra\cong \ZZ_p\times \ZZ_p$, $p\ne 3$ is a prime,
$a^\r=b$, $b^\r=a^{-1}b^{-1}$, $a^\s=b$, $b^\s=a$. In this case,
 $\M(G)\cong \ZZ_p$;

(iv) $G=V_1\times V_2 \times V_3$, $\a_i: V\to V_i$, $i=1,2,3$, are isomorphisms, $V$ is a finite
nonabelian simple group, $(v_1,v_2,v_3)^\s=(v_2,v_1,v_3)$, $(v_1,v_2,v_3)^\r=(v_3,v_1,v_2)$. In this
case, $\M(G)\cong V$;

(v) $G\cong P\Omega_8^+(q)$, $S$ is the group of graph automorphisms of $G$.
In this case, $\M(G)\cong M(q)$ is the simple Paige loop.
\end{lemma}

{\it Proof.} See \cite{dor,nv}.  $\blacktriangle$

Let $F$ be a field. An $FS$-module $V$ is called a {\it triality module} if $V$ is a group
with triality $S$. The representations $\varphi_i^{(\chi)}$ of $S$ corresponding to the
indecomposable $FS$-modules $V_i^{(\chi)}$ are shown in Table \ref{tab1}. They depend on
whether the characteristic $\chi$ of $F$ is $2$, $3$, or otherwise. We may assume $F$ to be a
prime field, since the decomposition field of $S$ is prime in any characteristic \cite{cr}.
The fact whether a module is a triality module is indicated by '$\checkmark$' in the last
column.

\begin{lemma} \label{indecomp}
Table \ref{tab1} holds.
\end{lemma}

\begin{table}[!htb]
\caption{Indecomposable $S_3$-modules \label{tab1}}
\vspace{5mm}
\begin{tabular}{cc@{}c@{}c@{}cc@{}c}
\hline
$\chi$ &\ \ $V^{\strut}$ \ \ & \ dim $V$ \ &\ \ \ $\varphi$ \ \ \
&$(12)\varphi$ & $(123)\varphi$ & Triality \\
\hline

$\ne 2,3$ & $V_1^{(0)\strut}$  & $1$ & $\varphi_1^{(0)}$& $1$ & $1$ & $\checkmark$ \\
          & $V_2^{(0)}$ & $1$ & $\varphi_2^{(0)}$ & $-1$ & $1$ & --- \\
 & $V_3^{(0)}$  &  $2$ & $\varphi_3^{(0)}$ &
$\left( \begin{array}{rr}
  0 & 1 \\  1 & 0
\end{array}  \right)$ &
$\left( \begin{array}{rr}
  0 & 1 \\ -1 & -1
\end{array}  \right)$ & $\checkmark$ \\[15pt]
\hline

$2$  & $V_1^{(2)\strut}$ &  $1$ & $\varphi_1^{(2)}$  & $1$ & $1$ & $\checkmark$ \\
  & $V_2^{(2)}$ &  $2$   & $\varphi_2^{(2)}$ &
$\left( \begin{array}{cc}
  0 & 1 \\  1 & 0
\end{array}  \right)$ &
$\left( \begin{array}{cc}
  0 & 1 \\  1 &  1
\end{array}  \right)$ & $\checkmark$ \\[15pt]
  & $V_3^{(2)}$  &  $2$  & $\varphi_3^{(2)}$ &
$\left( \begin{array}{cc}
  1 & 1 \\  1 & 0
\end{array}  \right)$ &
$\left( \begin{array}{cc}
  1 & 0 \\  0 & 1
\end{array}  \right)$ & --- \\[15pt]
\hline

$3$  & $V_1^{(3)\strut}$ &  $1$ & $\varphi_1^{(3)}$ & $1$ & $1$ & $\checkmark$ \\[5pt]
     & $V_2^{(3)}$ &  $1$ & $\varphi_2^{(3)}$ & $-1$ & $1$ & $\checkmark$ \\[5pt]
    & $V_3^{(3)}$ &  $2$ & $\varphi_3^{(3)}$  &
$\left( \begin{array}{rr}
  1 & 0 \\  0 & -1
\end{array}  \right)$ &
$\left( \begin{array}{rr}
  1 & 1 \\  0 &  1
\end{array}  \right)$ & $\checkmark$ \\[15pt]
 & $V_4^{(3)}$  &  $2$   & $\varphi_4^{(3)}$&
$\left( \begin{array}{rr}
  -1 & \ph{-}0 \\  0 & 1
\end{array}  \right)$ &
$\left( \begin{array}{rr}
  1 & 1 \\  0 &  1
\end{array}  \right)$ & $\checkmark$ \\[15pt]
 & $V_5^{(3)}$  & $3$  & $\varphi_5^{(3)}$ &
$\left( \begin{array}{rrr}
  1 & 0 & \ph{-}0 \\  0 & -1 & 0 \\ 0 & 0 & 1
\end{array}  \right)$ &
$\left( \begin{array}{rrr}
  1 & \ph{-}1 & 1 \\  0 & 1 & -1 \\ 0 & 0 & 1
\end{array}  \right)$ & $\checkmark$ \\[20pt]
  & $V_6^{(3)}$  & $3$   & $\varphi_6^{(3)}$&
$\left( \begin{array}{rrr}
  -1 & \ph{-}0 & 0 \\  0 & 1 & 0 \\ 0 & 0 & -1
\end{array}  \right)$ &
$\left( \begin{array}{rrr}
  1 & -1 & \ph{-}1 \\  0 & 1 & 1 \\ 0 & 0 & 1
\end{array}  \right)$ & --- \\[20pt]
\hline
\end{tabular}
\end{table}

{\it Proof.} The indecomposable $FS$-modules are well-known and can be readily determined (for
example, using \S64 in \cite{cr}). Namely, if $\chi\ne 2,3$ then every indecomposable
$FS$-module is irreducible by Maschke's theorem. If $\chi=3$ then there are exactly $6$
indecomposable $FS$-modules by Theorem (64.6) of \cite{cr} and the ones shown in Table
\ref{tab1} are easily seen to be indecomposable and pairwise nonisomorphic. If $\chi=2$ then
by the same theorem there are two non-faithful indecomposable modules: the trivial one and the reducible
(projective) one. The irreducible $2$-dimensional module is projective and is the only indecomposable
faithful $FS$-module in characteristic $2$, which follows from the fact that this module is the only
faithful component of the induced modules of the $2$-Sylow subgroup of $S$ (Theorem (63.8) of
\cite{cr}).

The triality of a module $V$ is equivalent to the condition that the element $(\s-1)(1+\r+\r^2)$ of $FS$
annihilates $V$, which is directly verified in each case. $\blacktriangle$

\begin{lemma} \label{syl_s_simple}
Let $G$ be an $S$-simple group with triality. Let $p$ be a Sylow prime for
$\M(G)$. Then $G$ has a $p$-Sylow $S$-subgroup.
\end{lemma}

{\it Proof.} We analyze the cases (i)--(v) of Lemma
\ref{s-sim}. If $G$ is as in cases (i)--(iii), the claim readily follows.
Let case (iv) hold. Then, for any $P\in Syl_p(V)$, the group $P_1\times P_2 \times P_3$ is the required
$p$-Sylow $S$-subgroup of $G$, where $P_i=\a_i(P)$, $i=1,2,3$. Suppose that $G\cong P\Omega_8^+(q)$ as
in case (v) and $S$ is the group of graph automorphisms of $G$. We will use the structure of some
$S$-subgroups of $G$ (for details, see
\cite{kle}). We have $|G|=\frac{1}{d^2}q^{12}(q^6-1)(q^4-1)^2(q^2-1)=n_1n_2n_3$, where the integers

$$ n_1=\frac{3}{(3,q)}q^{12}(q^2-1)^4, \quad n_2=
\frac{(3,q)(q^2\!+\!q\!+\!1)(q^2\!-\!q\!+\!1)}{3},
\quad n_3=\left(\frac{q^2+1}{d}\right)^2 $$
are pairwise coprime. By hypothesis, $p|n_1n_2$. If $p|n_2$ then $Syl_p(C_G(S))\subseteq Syl_p(G)$,
since $C_G(S)\cong G_2(q)$ and $n_2$ divides $|G_2(q)|=q^6(q^6-1)(q^2-1)$. Hence, any subgroup in
$Syl_p(C_G(S))$ is the required one.

Therefore, we may assume that $p|n_1$. If $p|q$ then consider an
$S$-invariant parabolic subgroup  $R_{s2}$  of $G$
(in the notation of \cite{kle}) of order
$\frac{1}{d^2}q^{12}(q-1)^4(q+1)^3$.
It has the following $S$-invariant
structure $R_{s2}=q^{9}\rtimes (SL_2(q)\circ SL_2(q)\circ SL_2(q)).q-1$ with
$S$ naturally permuting the three factors $SL_2(q)$, where $\circ$ denotes the central
product. The group $SL_2(q)\circ SL_2(q)\circ SL_2(q)$ has an obvious
$p$-Sylow $S$-subgroup whose preimage in $R_{s2}$ is the required subgroup.

If $p\nmid q$ and $p\ne 3$ then consider an $S$-subgroup $I_{+4}$ of order
$\frac{4}{d^2}q^{4}(q^2-1)^4$. It has the structure
$(SL_2(q)\circ SL_2(q)\circ SL_2(q)\circ SL_2(q)).d\rtimes K$, where
$K=\ZZ_2\times \ZZ_2$ and $K\rtimes S\cong S_4$ naturally
permutes the four factors $SL_2(q)$. It is easy to see that this group
contains a $p$-Sylow $S$-subgroup.

Finally, if  $p\nmid q$ and $p=3$ then let $q\equiv\ve (3)$, $\ve=\pm 1$. The required $3$-Sylow
subgroup is inside an $S$-subgroup $I_{\ve 2}$ of order $\frac{192}{d^2}(q-\ve)^4$. To see this,
consider an $\ve2$-decomposition of the Cayley algebra $\O(q)=V_1\oplus\ldots\oplus V_4$ which is also a
$\ZZ_2\times\ZZ_2$-grading with $\bm{1}\in V_1$ (see Section 4 in \cite{gz_max}). Let $R$ be the
centralizer of this decomposition in $\Omega_8^+(q)$ extended by the $4$-group $\la
(12)(34),(13)(24)\ra$ of permutations of the subspaces $\{V_i\}$, $i=1,\ldots,4$. This group has the
following structure $R=(\ZZ_{\frac{1}{d}(q-1)})^4.d^3.2^3.2^2$. In particular, the $3$-Sylow subgroup of
$R$ is characteristic. Take the group $T=\la (23),(34)\ra\cong S_3$ of permutations of $\{V_i\}$,
$i=2,3,4$ which acts identically on $V_1$ and induces automorphisms of $\O(q)$. Then the group $I_{\ve
2}=\ov{R}\rtimes \ov{T}$ is the image in $P\Omega_8^+(q)$ of $R \rtimes T\subseteq \Omega_8^+(q)$. It is
$S$-invariant with $S$ centralizing $\ov{T}$ (because $T$ is a group of automorphisms of $\O(q)$). It is
now clear that the $3$-Sylow subgroup $O_3(\ov{R})\rtimes \la\ov{(2,3,4)}\ra$ of $I_{\ve 2}$ is
$S$-invariant. $\blacktriangle$

\begin{lemma} \label{t1}
Let $G$ be a group with triality $S=\la \r,\s \ra$ and let $\M(G)=\{[g,\s]\ |\ \allowbreak g\in
G\}=(M,\,.\,)$ be the corresponding Moufang loop. Then, for every subloop $P\leqslant \M(G)$, the group
$Q=\la P\cup P^\r\ra$ is an $S$-subgroup of $G$ such that $\M(Q)=P$ and $[Q,S]=Q$.
\end{lemma}

{\it Proof.} By Theorem 1 in  \cite{gz_lagr}, $Q$ is an $S$-subgroup and $\M(Q)=P$. By Lemma \ref{dual},
$[Q,S]$ is the $S$-subgroup of $G$ generated by $P$. Hence, we have $Q\subseteq [Q,S]$. $\blacktriangle$

\begin{lemma} \label{ploop}
Let $M$ be a finite Moufang $p$-loop. Then

(i)\, $PsInn(M)$ is a $p$-group;

(ii)\, Any group with triality $G$ such that $\M(G)\cong M$
and $[G,S]=G$ is a $p$-group.
\end{lemma}

{\it Proof.} (i)\, The kernel of the natural epimorphism $\lm: PsInn(M)\rightarrow \I(M)$
which acts by $\lm: (A,a)\mapsto A$ is a subgroup of $Nuc(M)$ hence is a $p$-group.
The group $\I(M)$ is a subgroup of $Mlt(M)$ which is a $p$-group by Lemma VI.2.2 in
\cite{bru} and Theorem 4 in \cite{gla2}.

(ii)\, Let $G$ be as stated. Then $G$ is finite, since it is a quotient of the finite group $\D(M)$ (see
Corollary 3 in \cite{dor}). We have $G/Z_S(G)\cong \E(M)$ and $|\E(M)|=|PsInn(M)|\cdot |M|$ (see
\cite{gz_tri}) is a power of $p$ by (i). Hence, it remains to prove that $K=Z_S(G)$ is a $p$-group. Let
$P\in Syl_p(G)$.  We have $G=O_{p'}(K)\times P$, since $K$ is a central subgroup of $G$. The condition
$[G,S]=G$ now implies $O_{p'}(K)=1$. $\blacktriangle$

In the following lemma,
for an integer $n$ and a prime $p$, we denote by $n_p$ the
maximal power of $p$ dividing $n$.

\begin{lemma} \label{cont_syl}
Let $G$ be a group with triality $S=\la \r,\s\ra$, let $M=\M(G)$,
and let $p$ be a prime.

(i) If $N\leqslant G$ is an $S$-subgroup containing a $p$-Sylow
subgroup of $G$ then $|\M(N)|_p=|M|_p$.

(ii) If $P\in Syl_p(G)$ is $S$-invariant then $\M(P)\in Syl_p(M)$.
\end{lemma}
{\it Proof.} (i)\ Let $H=C_G(\s)$. We have $|\M(N)|=|N|/|N\cap H|$, and $|M|=|G|/|H|$. By
Lagrange's theorem \cite{gz_lagr}, $|M|/|\M(N)|=|G|/|NH|$ is an integer. However,
$|G|_p=|N|_p$ divides $|N|$ which, in turn, divides $|NH|$. Hence $|M:\M(N)|$ is coprime with
$p$.

(ii)\ Since $\M(P)$ is  a $p$-loop and $|M:\M(P)|$ is coprime with $p$ by (i), it follows that $\M(P)\in
Syl_p(M)$. $\blacktriangle$

As a consequence, we have an alternative proof independent of the classification
\cite{gz_max} of the following assertion (cf. Lemma \ref{syl_simple} above):

\begin{corollary} \label{alt}
If $p$ is a Sylow prime for $M(q)$ then $Syl_p(M(q))$ is nonempty.
\end{corollary}
{\it Proof.} We have  $M(q)=\M(G)$, where  $G=P\Omega_8^+(q)$ is $S$-simple. By Lemma
\ref{syl_s_simple}, $G$ contains a $p$-Sylow $S$-subgroup. The claim follows by (ii) of Lemma
\ref{cont_syl}. $\blacktriangle$

\begin{lemma} \label{help} Let $G=VP$ be a group with triality $S$
such that  $P$ is a $p$-group, $V \trianglelefteqslant G$, $p\nmid |V|$.
Suppose that $\M(\ov{P})$ is generated by at most $2$ elements,
where $\ov{P}=G/P$ satisfies $[\ov{P},S]=\ov{P}$.
Then $G$ has a $p$-Sylow $S$-subgroup.
\end{lemma}
{\it Proof.} Denote $M=\M(G)$. Let $\{m_i\}_{i\in I}$ generate $M$ modulo $W=\M(V)$. By hypothesis, we
may assume $|I| \leqslant 2$. In particular, $N=\big\la \{m_i\}_{i\in I}\big\ra$ is a group (since
Moufang loops are diassociative) and $M=WN$. Hence, any $R\in Syl_p(N)$ is in $Syl_p(M)$. By Lemma
\ref{t1}, $Q=\la R\cup R^\r \ra$ is an $S$-subgroup of $G$ such that $\M(Q)=R$ and $[Q,S]=Q$. Hence
$QV/V=\ov{P}$ in view of $[\ov{P},S]=\ov{P}$. By Lemma \ref{ploop}, $Q$ is the required $p$-Sylow
$S$-subgroup. $\blacktriangle$

\section{Proof of the theorem}

We are going to prove the following extended form of Theorem \ref{main}:
\begin{theorem} \label{ext}
Let $M$ be a finite Moufang loop and let $p$ be a prime.
Let $G$ be a group with triality such that $\M(G)=M$. Then the following
conditions are equivalent:

(i) $M$ has a $p$-Sylow subloop,

(ii) $p$ is a Sylow prime for $M$,

(iii) $G$ has an $S$-invariant $p$-subgroup $Q$ such that $\M(Q)\in Syl_p(M)$.
\end{theorem}

{\it Proof.} Let $k\geqslant 0$ be such that $p^k \|\, |M|$.

(i) $\Rightarrow$ (ii). Suppose that $M$ possesses a $p$-Sylow subloop $P$. We prove by
induction on the length of a composition series of $M$ that $p$ is a Sylow prime for $M$. If
$M$ is simple, this follows from Lemma \ref{syl_simple}. Let $N$ be a proper normal subloop of
$M$. It suffices to show that $p$ is a Sylow prime for both $N$ and $\ov{M}=M/N$. Note that
$P\cap N$ is a $p$-subloop of $N$ and $\ov{P}=\la P,N \ra/N$ is a $p$-subloop of $\ov{M}$ by

\be \label{third} \frac{\la P,N \ra}{N}\cong\frac{P}{P\cap N}, \ee
(see Theorem IV.1.5 in \cite{bru}). By Lagrange's theorem \cite{gz_lagr},
$|P|$ divides $|\la P,N \ra|$ which
in turn divides $|M|$. Hence (\ref{third}) also implies that the integers
$$\frac{|N|}{|P\cap N|}=\frac{|\la P,N \ra|}{|P|},\qquad
\frac{|\ov{M}|}{|\ov{P}|}=\frac{|M|}{|\la P,N \ra|}$$
are coprime with $p$. It follows that $P\cap N\in Syl_p(N)$
and $\ov{P}\in Syl_p(\ov{M})$.
By induction, $p$ is a Sylow prime for both $N$ and $\ov{M}$.

(ii) $\Rightarrow$ (iii). We now suppose that $p$ is a Sylow prime for $M=\M(G)$.
Proceed by induction on $|M|$.

If $Z_S(G)\ne 1$ then consider $G_0=G/Z_S(G)$. This group satisfies
$Z_S(G_0)=1$ and $\M(G_0)=M$. If we show that $P_0$ is an $S$-invariant
$p$-subgroup of $G_0$ such that $\M(P_0)\in Syl_p(M)$ then the preimage $P$
of $P_0$ in $G$ also satisfies $\M(P)\in Syl_p(M)$. If $P$ is not a $p$-group
then we take  $[P,S]$ which is  a $p$-group
by (ii) of Lemma \ref{ploop} and $\M([P,S])=\M(P)$.

Hence, we may assume that $Z_S(G)=1$. Clearly, we may also assume that $[G,S]=G$. Take a
minimal normal $S$-subgroup $V$ of $G$. Then $V$ is characteristically simple (see Theorem 1.5
in Chap. 2 of \cite{gor}). Two cases are possible:

(1) $p\mid |V|$. Then we may assume that $V$ is nonabelian (otherwise $V$ is a $p$-group and we apply
induction for $G/V$, which is possible, since $p$ is a Sylow prime for $\M(G/V)$ and $|\M(G/V)|<|M|$).
We have $V=V_1\times\ldots\times V_s$ is the product of isomorphic nonabelian simple groups $V_i$'s.
Since $p$ is Sylow for $\M(V)$ and $V$ is a direct product of $S$-simple groups, it follows by Lemma
\ref{syl_s_simple} that $V$ has a $p$-Sylow $S$-subgroup $W$.

By assumption $W\ne 1$. We have  $G=VN$, where $N=N_G(W)$ is $S$-invariant. Observe that
$N$ contains a $p$-Sylow subgroup of $G$. Moreover $N<G$, since $W\ntrianglelefteqslant
V$;  and $\M(N)<M$, since $[G,S]=G$. Also note that $p$ is Sylow for $\M(N)$, since
$W\trianglelefteqslant N\cap V \trianglelefteqslant N$ is an $S$-invariant series for
$N$ and $W$ is a $p$-group, $(N\cap V)/W$ is a $p'$-group, and $N/(N\cap V)\cong G/V$.
By induction, there  exists a $p$-subgroup $P$ in $N$ such that $\M(P)\in Syl_p(\M(N))$.
By (i) of Lemma \ref{cont_syl}, $\M(P)\in Syl_p(M)$.

(2) $p\nmid |V|$. By induction, we may assume that $\ov{P}=G/V$ is a $p$-group. It is sufficient to show
that $G$ has a nontrivial $S$-invariant $p$-subgroup $P_0$ such that $\ov{P_0}=P_0V/V$ is normal in
$\ov{P}$. Indeed, if such a subgroup exists then the condition $Z_S(G)=1$ implies that either
$\M(P_0)\ne 1$ or $P_0 \ntrianglelefteqslant G$. In both cases we can use induction for $N/P_0$, where
$N=N_G(P_0)$. This is because $|\M(N/P_0)|<|M|$ and $p$ is Sylow for $\M(N/P_0)$, since $N/P_0$ is an
extension of the $p'$-group $C_V(P_0)$ by the $p$-group $\ov{P}/\ov{P_0}$. By induction, $N/P_0$
contains an $S$-invariant $p$-subgroup $\ov{P_1}$ such that $\M(\ov{P_1})\in Syl_p(\M(N/P_0))$. Then the
full preimage $P_1$ of $\ov{P_1}$ in $N$ in the required $S$-subgroup, since $\M(P_1)\in
Syl_p(\M(N))\subseteq Syl_p(M)$.

Let $\ov{Z}=Z(\ov{P})$ and suppose that $\M(\ov{Z})\ne 1$.
The existence of $P_0$ is this case is easy.
Take $z\in \M(G)$ such that $1\ne \ov{z}\in \M(\ov{Z})$, where $\ov{z}=Vz$.
The group $G_0=\la V,z,z^\r \ra$ satisfies the conditions of Lemma \ref{help}.
We take $P_0$ to be a $p$-Sylow $S$-subgroup of $G_0$.
Then $\ov{P_0}\trianglelefteqslant\ov{P}$ as a central
subgroup of $\ov{P}$.

Suppose that $\M(\ov{Z})=1$. Then $\M(\ov{Z_1})\ne 1$, where
$\ov{Z_1}/\ov{Z}=Z(\ov{P}/\ov{Z})$, since otherwise we would have
$\ov{Z}_1\leqslant Z_S(\ov{P})\leqslant \ov{Z}$ (the latter inclusion follows
from $[\ov{P},S]=\ov{P}$ and (ii) of Lemma \ref{dual}).
Take $a\in \M(G)$ such that $1\ne \ov{a}\in \M(\ov{Z_1})$
and put $A=\la a,a^\r \ra$. Then $A$ is an  $S$-subgroup of $G$
not contained in $Z$, where $Z$ is the full preimage of $\ov{Z}$ in $G$.

The elementary abelian group $\ov{P}/\Phi(\ov{P})$ is the direct product
$U_1\times \ldots\times U_t$
of indecomposable triality $\F_p S$-modules $U_i$'s.
The condition $[\ov{P},S]=\ov{P}$ implies that the $U_i$'s
are at most $2$-dimensional and, depending on $p$, are isomorphic
to one of the modules $V_3^{(0)}$, $V_2^{(2)}$, $V_2^{(3)}$, $V_4^{(3)}$
from Table \ref{tab1}. Moreover, we have
$\M(U_i)=\la u_i\ra$ is cyclic of order $p$ and
$U_i=\la u_i,u_i^\r \ra$  (observe that
$u_i^\r=u_i$ if $U_i \cong V_2^{(3)}$).
Let $w_i$ be corresponding preimages
of $u_i$ in $\M(G)$. Then $G$ is generated modulo $V$
by the $S$-subgroups $W_i=\la w_i,w_i^\r \ra$. Since
$Z_1\geqslant A\nleqslant Z$,
where $Z_1$ is the full preimage of $\ov{Z_1}$ in $G$, it follows that
there exists $i_0$ such that $W=[A,W_{i_0}]\nleqslant V$.
On the other hand, $W\leqslant Z$, since $A\leqslant Z_1$.
Denote  $G_0=\la V,A,W_{i_0} \ra$ and  $\ov{G_0}=G_0/V$.
Then the $S$-subgroup $G_0$ satisfies
the conditions of Lemma \ref{help}, since the images of $a$
and $w_{i_0}$ in $\ov{G_0}$
generate $\M(\ov{G_0})$ as a loop and $\ov{G_0}$
as an $S$-group (whence the condition
$[\ov{G_0},S]=\ov{G_0}$).
Let $P_1$ be a $p$-Sylow $S$-subgroup of $G_0$.
Put $P_0=P_1\cap Z$. It remains to observe that $P_0\ne 1$, since
$P_0\cap WV=P_1 \cap WV \ne 1$ in view of  $W\nleqslant V$.

(iii) $\Rightarrow$ (i). Obvious. $\blacktriangle$

We note that it is not true in general that if $G$ is a group with triality and $p$ is Sylow for $\M(G)$
then $G$ contains a $p$-Sylow $S$-subgroup. For example, let $G=S$ on which $S$ acts by inner
automorphisms. Then $2$ is a Sylow prime for $\M(G)\cong \ZZ_3$, but $G$ does not have a $2$-Sylow
$S$-subgroup.

\section{The group-type radical of Moufang loops}

A finite Moufang loop $M$ is said to be a loop of {\it group type} if all composition factors
of $M$ are groups. For example, all solvable Moufang loops are loops of group type. It is
clear that the normal subloop of $M$ generated by two normal subloops of group type is again a
loop of group type. Hence we have
\begin{pro} Every finite Moufang loop has a unique maximal normal sub\-loop of group type.
\end{pro}
We denote this maximal normal subloop of group type by $Gr(M)$. It is obvious that $Gr(M/Gr(M))=1$,
hence we call $Gr(M)$ the {\it group-type radical} of $M$.

For $q$ odd, the simple Paige loop $M(q)$ has a two-fold extension isomorphic to the
loop $PGL(\O(q))$, where $\O(q)$ is the Cayley algebra over $\F_q$ (see Section 4 in
\cite{gz_max}). We denote this extension by $M(q).2$. Also define

$$ \widehat{M(q)}=\left\{\begin{array}{ll}
M(q).2, & \text{if}\ q\  \text{is odd}\\
M(q), & \text{if}\ q\ \text{is even}
\end{array}
  \right.$$
It can be seen that the group $InnDiag(P\Omega_8^+(q))$ of inner-diagonal automorphisms of
$P\Omega_8^+(q)$ is a group with triality $S$ corresponding to $\widehat{M(q)}$, where $S$ is
the group of graph automorphisms of $P\Omega_8^+(q)$. It is known that the factor group
$$InnDiag(P\Omega_8^+(q))/P\Omega_8^+(q)$$ is trivial for $q$ even and isomorphic to $\ZZ_2\times
\ZZ_2$ for $q$ odd. Moreover, in the latter case, $S$ acts
nontrivially on $\ZZ_2\times \ZZ_2$ so that $\M(\ZZ_2\times \ZZ_2)=\ZZ_2$, whence
$|\widehat{M(q)}|=2|M(q)|$.

The following assertion is the main result of this section (cf. Theorem \ref{thb} in the
introduction):
\begin{theorem} \label{MT}Let $M$ be a finite Moufang loop and let $N$ be a minimal normal
subloop of $M$. Then we have

(i)\ $N$ is the direct product of isomorphic simple loops. Furthermore, if $N$ is nonassociative then it
is simple.

(ii)\ Suppose, additionally, that $M$ satisfies $Gr(M)=1$. Then by (i) every minimal
normal subloop of $M$ is a simple Paige loop $M(q)$ for some $q$. Denote $Soc(M)=\prod
M(q)$ and $\widehat{Soc(M)}=\prod \widehat{M(q)}$, where both products are taken over
all minimal normal subloops of $M$. Then
$$Soc(M)\trianglelefteqslant M \trianglelefteqslant \widehat{Soc(M)}.$$
In particular, $M/Soc(M)$ is an elementary abelian $2$-group.

\end{theorem}
 {\bf Proof.} Suppose that $M$  and $N$ are as stated.

(i)\ Let $G$ be a group with triality $S=\la \r,\s\ra$ such that $\M(G)=M=\{x^{-1}x^{\s}\ |\ x\in G\}$,
$[G,S]=G$, and let $Q_0$ be a minimal normal $S$-subgroup of $G$ corresponding to $N$. Since $Q_0$ is
characteristically simple, we have $Q_0=Q_1\times \ldots \times Q_n$, where $Q_i, 1\leqslant i\leqslant
n$, are isomorphic simple groups. Note that we also have the decomposition $Q_0=R_1\times \ldots
\times R_k$, where $R_j$, $1\leqslant j\leqslant k$, are $S$-simple groups if $Q_0$ is
nonabelian or indecomposable $S$-modules if $Q_0$ is abelian; moreover, $N=\M(Q_0)=\M(R_1)\times \ldots
\times \M(R_k)$. If $Q_1\not\cong P\Omega_8^+(q)$ for any $q$ then Lemmas \ref{indecomp} and \ref{s-sim}
imply that all $\M(R_j)$'s are either trivial or isomorphic simple groups. Hence, we may
assume that $Q_1\cong P\Omega_8^+(q)$, with $q=p^m$, $p$ a prime. We denote $I=\{i\ |\
Q_i^{S}=Q_i\}$. If either $I=\varnothing$ or $[Q_i,S]=1$ for all $i\in I$ then each $R_j$ is
either an $S$-group from item (iv) of Lemma \ref{s-sim} with $V_t= P\Omega_8^+(q)$, $t=1,2,3$,
or is isomorphic to $P\Omega_8^+(q)$ with trivial $S$-action. Hence, in this case, $N$ is the
direct product of several $P\Omega_8^+(q)$'s and the claim holds. Consequently, we may assume
that $1\in I$ and $\M(Q_1)=M(q)\ne 1$. We prove that $n=1$ in this case. Assume the contrary.
Denote $K=N_G(Q_1)$. Observe that $K$ is $S$-invariant. We have $K\ne G$, since otherwise
$Q_1$ would be a normal $S$-subgroup contrary to the choice of $Q_0$. Note that we cannot have
the inclusion $M\subseteq K$, since $G$ is $S$-generated by $M$ (i.e.,  $G=\la\cup_{\tau\in S}
M^\tau \ra$, which follows from $[G,S]=G$ by (iii) of Lemma \ref{dual}). Hence, there exists
$x\in M\setminus K$ such that $x^l\in K$ for a prime $l$. Let $y=x^{\rho}$. Then we have

\begin{equation}\label{trm}
x^{\s}=x^{-1},\,y^{\rho}=x^{-1}y^{-1},\,y^{\s}=xy,\,[x,y]=1\end{equation}

in view of (\ref{m-1}) and (iv) of Lemma \ref{dual}. Denote $Q_{(i,j)}=Q_1^{x^iy^j},$ $i,j\in
{\ZZ}/l{\ZZ},$ $Q=Q_{(0,0)}=Q_1$. Suppose that $x^ay^b\in K,$ then by (\ref{trm}) we have
$(x^ay^b)^{\s}=x^{b-a}y^b\in K$ and $(x^ay^b)^{\rho}=x^{-b}y^{a-b}\in K$. Hence $3a\equiv 3b\equiv 0
\pmod{l}$.

I. $l>3.$ Then $a\equiv b\equiv 0 \pmod{l}$ and all groups $Q_{(i,j)}$ are different. By
(\ref{trm}) we obtain

\begin{equation}\label{sigma}
Q_{(i,j)}^{\s}=Q_{(j-i,j)},\,Q_{(i,j)}^{\rho}=Q_{(-j,i-j)}.
\end{equation}
Since $Q_0$ is a group with triality, the orbits of $S$ on the set $\{Q_1, ... ,Q_n\}$ are all of
size $1$ or $3$. Hence by (\ref{sigma}) the set
$$\{(i,j),(j-i,j),(-j,i-j),(i,i-j),(-j,-i),(j-i,-i)\}$$
has $3$ elements or less for all $i,j\in {\ZZ}/l{\ZZ}$. However, this is possible only in the case
$l=2.$

II. $l=3.$ Since $x^3,y^3\in K$ we have by (\ref{trm})

$(Q^{xy^{-1}})^{\rho}=Q^{xy^2}=(Q^{y^3})^{xy^{-1}}=Q^{xy^{-1}},$

$(Q^{xy^{-1}})^{\s}=Q^{x^{-2}y^{-1}}=(Q^{x^{-3}})^{xy^{-1}}=Q^{xy^{-1}}.$

Similarly, we can prove that the set $X=\{Q^x,Q^{x^2},Q^y,Q^{y^2},Q^{xy},Q^{x^2y^2}\}$ is an
$S$-orbit. Since $|X|>1$ and $Q$ is a simple group, we have $|X|=3$. If $Q^x=Q^{y^2}$, an
application of $\r$ gives $Q^y=Q^{xy}$, which implies $Q^x=Q$ in contradiction with the choose of
$x$. The cases $Q^x=Q^{x^2}$ or $Q^x=Q^{xy}$ are treated similarly. In the other cases: $Q^x=Q^y$,
$Q^x=Q^{x^2y^2}$ we have $|X|=2$, a contradiction.

III. $l=2.$ Then the subgroup $Q\times
Q^x\times Q^y\times Q^{xy}$ is $S$-invariant and

$(Q^x)^{\s}=Q^x,\,(Q^y)^{\s}=Q^{xy},\,(Q^x)^{\rho}=Q^{y},\,
(Q^{y})^{\rho}=Q^{xy},\,(Q^{xy})^{\rho}=Q^x.$

Let $a\in Q$. We apply the identity of triality (\ref{trial}) to $g=ax$. We have

$r=g^{-1}g^{\s}=x^{-1} (a^{-1}a^{\s}){x^{-1}}$ and

$1=rr^\r r^{\r^2}= x^{-1}(a^{-1}a^{\s}) x^{-1}x^{-\r}(a^{-1}a^{\s})^\r x^{-\r}
x^{-\r^2}(a^{-1}a^{\s})^{\r^2} x^{-\r^2}=$

$x^{-1}(a^{-1}a^{\s}) x^{-1}y^{-1}(a^{-1}a^{\s})^\r x (a^{-1}a^{\s})^{\r^2} xy=$

$x^{-1}(a^{-1}a^{\s}) x^{-1}y^{-1}(a^{-1}a^{\s})^\r x^2y \big{(} (a^{-1}a^{\s})^{\r^2} \big{)}^ {xy}=$

$x^{-1}(a^{-1}a^{\s}) x \big{(} (a^{-1}a^{\s})^\r \big{)}^{x^2y}  \big{(} (a^{-1}a^{\s})^{\r^2} \big{)}^ {xy}=$

$\big{(}a^{-1}a^{\s}\big{)}^x \big{(} (a^{-1}a^{\s})^\r \big{)}^{x^2y} \big{(}
(a^{-1}a^{\s})^{\r^2} \big{)}^ {xy}$.

Since $x^2\in K$, the last expression is in $Q^x\times Q^y \times Q^{xy}$. Hence,
$a^{-1}a^{\s}=1$. By arbitrariness of $a$, we have $\M(Q)=\{a^{-1}a^{\s} | a\in Q\}=1$, a
contradiction.  This proves item (i).

(ii) Suppose that $Gr(M)=1$. Let $G$ be as above with the additional condition that $Z_S(G)=1$.
Let $T$ be a minimal normal
$S$-subgroup of $G$ corresponding to $Soc(M)$. By the above, $T=\prod_q P\Omega_8^+(q)$ with $S$
acting on each factor of $T$ by graph automorphisms. Observe that $C_G(T)$ is a normal $S$-subgroup
of $G$. Denote $M_0=\M(C_G(T))\trianglelefteqslant M$. Since $C_G(T)\cap T=1$, we have $M_0\cap
Soc(M)=1$, which implies $M_0=1$ and, therefore, $C_G(T)=1$ in view of $Z_S(G)=1$. It follows that
$T\trianglelefteqslant G\leqslant Aut(T)$. We have
$$Aut(T)=\prod_i
\big{(}\underbrace{Aut(P\Omega_8^+(q_i))\times\ldots
             \times Aut(P\Omega_8^+(q_i))
            }_{n_i}
\big{)}\rtimes S_{n_i},
$$
where $S_{n_i}$ is the symmetric group of degree $n_i$ naturally permuting the isomorphic factors
$Aut(P\Omega_8^+(q_i))$'s. Since $S$ acts on each $P\Omega_8^+(q_i)$ in the same way, this action
commutes with the action of $S_{n_i}$. Moreover, $S$ commutes with the field automorphisms of
$P\Omega_8^+(q_i)$ (for details, see p. 181 in \cite{kle}, where the structure of $Aut(P\Omega_8^+(q))$
is discussed). Thus, the condition $[G,S]=G$ forces $G$ to be an $S$-subgroup of $\prod_j
D_j\rtimes\Gamma_j\leqslant Aut(T)$, where $D_j=InnDiag(P\Omega_8^+(q_j))$ and $\Gamma_j$ is the group
of graph automorphisms of $P\Omega_8^+(q_i)$. Since the $j$th projection $G\to D_j\rtimes\Gamma_j$
commutes with the action of $S$, its image $G_j$ must be a subgroup with triality in
$D_j\rtimes\Gamma_j$. However $\M(G_j)$, containing $M(q)$ as a subloop, does not have exponent $3$. By
Lemma \ref{tripl}, we must have $G_j\leqslant D_j$. Hence, $G_j$ is either $P\Omega_8^+(q_j)$ or
$InnDiag(P\Omega_8^+(q_j))$, and $\M(G_j)$ is either $M(q)$ or $\widehat{M(q)}$, accordingly.
$\blacktriangle$

Observe that part III of the proof of item (i) implies, in particular,
the following useful fact:

\begin{pro} Let the symmetric group $S_4$ act on the direct
product of isomorphic groups $G=G_1\times G_2\times G_3\times G_4$ in such a way that
$(G_i)^\tau=G_{i^\tau}$ for all $i=1,\ldots,4$ and $\tau \in S$. Let $S_4=S\ltimes N$, where
$S=\la \s=(12),\r=(123)\ra$ and $N=\la (12)(34), (13)(24) \ra$. If the semidirect product
$N\ltimes G$ is a group with triality $S$ then $G_4\leqslant C_G(S)$.
\end{pro}

The following lemma will be used in the next section:

\begin{lemma} \label{gr_prop} Let $M$ and $N$ be Moufang loops.

(i)\ If $\varphi: M\to N$ is a homomorphism then
$Gr(M)^\varphi\leqslant Gr(M^\varphi)$, where the inclusion can be proper;

(ii)\ If $Gr(M)=1$ and $N\trianglelefteqslant M$ then $Gr(N)=1$.
\end{lemma}
{\it Proof.} The inclusion in (i) readily follows from the definition. If we take $M=\widehat{M(3)}$ and
$N=\mathbb{Z}_2$ then $N$ is a homomorphic image of $M$; however, $Gr(M)=1$ and $Gr(N)=N$. To show (ii),
define $\pi_i$ to be the projection of $\widehat{Soc(M)}$ to its $i$th direct factor isomorphic to
$\widehat{M(q)}$ (see item (ii) of Theorem \ref{MT}). Then $\pi_i$ maps $N$ to a normal subloop of
$\pi_i(M)$ which is either $M(q)$ or $\widehat{M(q)}$. In any case, $Gr(N^{\pi_i})=1$. By (i), this
implies $Gr(N)^{\pi_i}=1$ for all $i$ and the claim follows. $\blacktriangle$

\section{$p$-Subloops for non-Sylow primes $p$}

The main theorem about the existence of $p$-Sylow subloops suggests the following natural question: How
large can a $p$-subloop of a Moufang loop $M$ be if $p$ is not a Sylow prime for $M$? In this section,
we show that it can be as large as possible. To give a more precise statement, we need another
definition.

Suppose that $M$ is a finite Moufang loop and $M_i$, $1\leqslant i \leqslant l$, are the
composition factors of $M$. Let $p$ be any prime. Take $P_i\in Syl_p(M_i)$ if $p$ is a Sylow
prime for $M_i$ and put $P_i=1$, otherwise. Then a {\it quasi-$p$-Sylow} subloop of $M$ is a
subloop of order $\prod_{1\leqslant i\leqslant l} |P_i|$. By Lemma~\ref{syl_simple}, a
quasi-$p$-Sylow subloop, if exists, must be a maximal $p$-subloop of $M$ and its order must be
the maximal order of all $p$-subloops of $M$. Denote by $Q\mbox{-}Syl_p(M)$ the set of all
quasi-$p$-Sylow subloops of $M$. Clearly, if $p$ is a Sylow prime for $M$ then
$Q\mbox{-}Syl_p(M)=Syl_p(M)$. Our main result of this section is the following assertion (cf.
Theorem \ref{tha} in the introduction):

\begin{theorem} \label{quasi} Let $p$ be a prime and $M$ a finite Moufang loop.
Then $M$ contains a quasi-$p$-Sylow subloop.
\end{theorem}
{\it Proof.} We may assume that $p$ is non-Sylow for $M$, since otherwise the result follows from
Theorem \ref{ext}. This implies that $p>3$. We proceed by induction on the composition length $l$ of
$M$. If $l=1$ then the claim holds by definition. Assume that $l>1$. Let $N$ be a minimal normal subloop
of $M$. By induction, there exists $\ov{P}\in Q\mbox{-}Syl_p(\ov{M})$, where $\ov{M}=M/N$. Let $P$ be
the full preimage of $\ov{P}$ in $M$. If $p$ is a Sylow prime for $N$ then it is such for $P$ as well.
In this case, $P$ contains a $p$-Sylow subloop by Theorem \ref{ext} which is obviously quasi-$p$-Sylow
for $M$. Hence, we may assume that $p$ is non-Sylow for $N$. But then $N$ must be nonassociative and, by
Theorem
\ref{MT}, $N$ is simple. Note that a quasi-$p$-Sylow subloop of $P$, if exists, must have
order $|\ov{P}|$. Let $R=Gr(P)$. Since $p$ is odd, it is easy to conclude by (ii) of Theorem \ref{ext}
that $P/R\cong N$ and $|R|=|\ov{P}|$. Hence, $R\in Q\mbox{-}Syl_p(P)\subseteq Q\mbox{-}Syl_p(N)$.
$\blacktriangle$

A corresponding fact about groups with triality can be stated as follows:

\begin{corollary} Let $p$ be a prime and let $G$ be a group with triality $S$. Then
$G$ possesses an $S$-invariant $p$-subgroup $P$ such that $\M(P)\in Q\mbox{-}Syl_p(\M(G))$.
\end{corollary}

{\it Proof.} This assertion is a generalization of the implication (i) $\Rightarrow$ (iii) of Theorem
\ref{ext} and can be proved as follows. Denote $M=\M(G)$. We identify $M$ with the subset $\{ [x,\s]\ |\
x\in G\}\subseteq G$. Let $N\in Q\mbox{-}Syl_p(M)$, which exists by Theorem \ref{quasi}. By Lemma
\ref{t1}, $P=\la N\cup N^\r \ra$ is an $S$-subgroup of $G$ such that $\M(P)=N$ and $[P,S]=P$. By (ii) of
Lemma \ref{ploop}, $P$ is the required $p$-subgroup. $\blacktriangle$

Let $M$ be a Moufang loop and let $p$ be a prime. Denote by $Gr_p(M)$ the product of all
normal subloops of $M$ for which $p$ is a Sylow prime. The properties of $Gr_p(M)$
are as follows:

\begin{pro} \label{gr_p_prop} (i)\ $Gr_p(M)\trianglelefteqslant M$ and $p$ is a Sylow
prime for $Gr_p(M)$;

(ii)\ The factor loop $M/Gr_p(M)$ contains no elements of order $p$;

(iii)\ All $p$-subloops of $M$ are contained in $Gr_p(M)$;

(iv)\ $Syl_p(Gr_p(M))=Q\mbox{-}Syl_p(M)$;

(v)\ $Gr(M)=Gr(Gr_p(M))=\bigcap_p Gr_p(M)$.

\end{pro}
{\it Proof.} Item (i) follows directly from the definition. We may henceforth assume that $p$ is
non-Sylow for $M$, since otherwise $Gr_p(M)=M$ and the claim trivially holds. In particular, we have
$p>3$. Observe that $Gr(M)\leqslant Gr_p(M)$. Hence, $M/Gr_p(M)\cong \overline{M}/Gr_p(\overline{M})$,
where $\overline{M}=M/Gr(M)$, and it is easy to conclude by Theorem \ref{MT} that $Soc(M/Gr_p(M))$ must
be the product of $M(q)$ for which $p$ is non-Sylow. Hence, $M/Gr_p(M)$ does not contain elements of
order $p$ by Lemma \ref{syl_simple} and (ii) follows. If $P$ is a $p$-subloop of $M$ then
$PGr_p(M)/Gr_p(M)$ is a $p$-subloop of $M/Gr_p(M)$ and hence must be trivial by (ii), which implies
(iii). The inclusion $Syl_p(Gr_p(M))\subseteq Q\mbox{-}Syl_p(M)$ holds by the definition of a
quasi-$p$-Sylow subloop, since the composition factors of $M/Gr_p(M)$ are $M(q)$ for which $p$ is
non-Sylow and, possibly, cyclic groups of order $2$, as we just explained in proving (ii). The reverse
inclusion holds by (iii) and the fact that quasi-$p$-Sylow subloops are $p$-subloops of maximal order.
Finally, we show that (v) holds. The inclusions $Gr(M)\leqslant Gr(Gr_p(M))$ and $Gr(M) \leqslant
\bigcap_p Gr_p(M)$ are obvious from $Gr(M)\leqslant Gr_p(M)$. Since $Gr_p(M)/Gr(M)$ is a normal subloop
of $M/Gr(M)$, its group-type radical must be trivial by Lemma \ref{gr_prop}. It follows easily that
$Gr(Gr_p(M))\leqslant Gr(M)$. Now denote by $R$ the image of $\bigcap_p Gr_p(M)$ in $M/Gr(M)$. Since $R$
is a normal subloop, it follows by (ii) of Lemma \ref{gr_prop} that every minimal normal subloop of $R$
is $M(q)$ for some $q$. However, for every $M(q)$, there is a non-Sylow prime. Hence $R$ has no
nontrivial normal subloops and thus must be trivial. $\blacktriangle$

This proposition shows that $Gr_p(M)$ can be viewed as a $p$-analog of the group-type
radical $Gr(M)$. In particular, the study of embeddings of $p$-subloops into each other
and determining the number of quasi-$p$-Sylow subloops of $M$ can be reduced to the
case when $p$ is a Sylow prime for $M$.

\end{document}